\documentclass[11pt,a4paper]{article}
\usepackage[top=1in, bottom=1.25in, left=1in, right=1in]{geometry}
\usepackage[latin1]{inputenc}
\usepackage[english]{babel}
\usepackage{graphicx,epstopdf}
\usepackage{caption, subcaption,float}
\usepackage{amsthm}
\usepackage{color}
\usepackage{tikz}
\usepackage{relsize}
\usepackage[ruled, linesnumberedhidden,vlined]{algorithm2e} 
\usepackage[round,sort,comma]{natbib}
\usepackage{hyperref}
\hypersetup{
    colorlinks=true,       
    linkcolor=blue,          
    citecolor=olive,       
    filecolor=magenta,      
    urlcolor=blue          
}

\usetikzlibrary{matrix}
\usepackage{amsmath,amssymb,dsfont}
\begin{document}

\title{On the minimal positive standardizer of a parabolic subgroup of an Artin-Tits group}
\date{\today }
\author{Mar\'{i}a Cumplido\\  Universit\'{e} Rennes 1, Universidad de Sevilla}

\maketitle
\newtheorem{theorem}{Theorem}
\newtheorem*{theorem*}{Theorem 4}
\newtheorem{lemma}[theorem]{Lemma}
\newtheorem{proposition}[theorem]{Proposition}
\newtheorem{corollary}[theorem]{Corollary}

\theoremstyle{definition}
\newtheorem{definition}[theorem]{Definition}

\theoremstyle{remark}
\newtheorem{remark}[theorem]{Remark}

\newcommand{\St}[1]{\operatorname{St}(\mathcal{#1})}
\newcommand{\po}{\preccurlyeq}
\newcommand{\wedger}{\wedge^\Lsh}
\newcommand{\veer}{\vee^\Lsh}
\newcommand{\bigveer}[1]{\underset{#1\,\,\,}{{\bigvee}^{\mathlarger{\Lsh}}}}
\newcommand{\alphae}{\alpha_{\text{ext}}}
\newcommand{\gammae}{\gamma_{\text{ext}}}
\newcommand{\notpo}{\npreceq}
\newcommand{\so}{\succcurlyeq}
\newcommand{\notso}{\nsucceq}

\begin{abstract}
The minimal standardizer of a curve system on a punctured disk is the minimal positive braid that transforms it into a system formed only by round curves. We give an algorithm to compute it in a geometrical way. Then, we generalize this problem algebraically to parabolic subgroups of Artin-Tits groups of spherical type and we show that, to compute the minimal standardizer of a parabolic subgroup, it suffices to compute the $pn$-normal form of a particular central element. 
\end{abstract}

\section{Introduction}

\noindent Let $D$ be the disk in $\mathbb{C}$ with diameter the real segment $[0,n+1]$ and let $D_n= D\setminus \{1,\dots, n\}$ be the $n$-punctured disk. The $n$-strand braid group, $B_n$, can be identified with the mapping class group of $D_n$ relative to $\partial D_n$. $B_n$ acts on the right on the set of isotopy classes of simple closed curves in the interior of $D_n$. The result of the action of a braid $\alpha$ on the isotopy class~$\mathcal{C}$ of a curve~$C$ will be denoted by~$\mathcal{C}^{\alpha}$ and it is represented by the image of the curve $C$ under any automorphism of~$D_n$ representing $\alpha$. We say that a curve is \emph{non-degenerate} if it is not homotopic to a puncture, to a point or to the boundary of $D_n$, in other words, if it encloses more than one and less than $n$ punctures. A \emph{curve system} is a collection of isotopy classes of disjoint non-degenerate simple closed curves, pairwise non-isotopic.

\medskip

Curve systems are very important as they allow to use geometric tools to study braids. From Nielsen-Thurston theory \citep{Thurston1988}, every braid can be decomposed along a curve system, so that each component becomes either periodic or pseudo-Anosov. The simplest possible scenario appears when the curve is \emph{standard}:

\begin{definition}A simple closed curve in $D_n$ is called \textbf{standard} if it is isotopic to a circle centered at the real axis. A curve system containing only isotopy classes of standard curves is called standard. \end{definition}

Every curve system can be transformed into a standard one by the action of a braid, as we shall see. Let $B_n^{+}$ be the submonoid of $B_n$ of positive braids, generated by $\sigma_1, \ldots , \sigma_{n-1}$ \citep{Artin1947}. We can define a partial order $\po$ on $B_n$, called \emph{prefix order}, as follows: for $\alpha,\beta\in B_n$, $\alpha \po \beta$  if there is $\gamma \in B_n^{+}$ such that $ \alpha \gamma = \beta$. This partial order endows $B_n$ with a lattice structure, i.e., for each pair $\alpha, \beta \in B_n$, their gcd $\alpha\, \wedge \beta$ and their lcm $\alpha\, \vee \beta$ with respect to $\po$ exist and are unique. Symmetrically, we can define the \emph{suffix order} $\so$ as follows: for $\alpha,\beta\in B_n$, $\beta \so \alpha$  if there is $\gamma \in B_n^{+}$ such that $ \gamma\alpha = \beta$. We will focus on $B_n$ as a lattice with respect to $\po$, and we remark that $B_n^+$ is a sublattice of $B_n$. In 2008, Lee and Lee proved the following:

\begin{theorem}[\citealp{Lee2008}]\label{Lee}
Given a curve system $\mathcal{S}$ in $D_n$,  its standardizer $$\St{S}=\{\alpha\in B_n^{+}\, :\, \mathcal{S}^{\alpha}\text{  is standard }\}$$ is a sublattice of $B_n^{+}$. Therefore, $\St{S}$ contains a unique $\po$-minimal element. 
\end{theorem}

The first aim of this paper is to give a direct algorithm to compute the $\po$-minimal element of $\St{S}$, for a curve system $\mathcal{S}$. The algorithm, explained in Section~\ref{salgo}, is inspired by Dynnikov and Wiest algorithm to compute a braid given its curve diagram \citep{Dynnikov2007} and the modifications made in \citep{Caruso2013}. 

\bigskip

The second aim of the paper is to solve the analogous problem for Artin-Tits groups of spherical type. 

\begin{definition}
Let $S$ be a finite set and $M=(m_{i,j})_{i,j\in S}$ a symmetric matrix with $m_{i,i}=1$ and $m_{i,j}\in\{2,\dots, \infty \}$ for $i\neq j$. Let $\Sigma=\{\sigma_i\,|\, i\in S\}$. The Artin-Tits system associated to $M$ is $(A,\Sigma)$, where $A$ is a group (called Artin-Tits group) with the following presentation 
$$A=\langle \Sigma \,|\, \underbrace{\sigma_i\sigma_j\sigma_i\dots}_{m_{i,j} \text{ elements}}=\underbrace{\sigma_j\sigma_i\sigma_j\dots}_{m_{i,j} \text{ elements}} \forall i,j\in S,\, i\neq j,\, m_{i,j}\neq \infty \rangle.$$
\end{definition}
\noindent For instance, $B_n$ has the following presentation \citep{Artin1947}
$$B_n=\left\langle \sigma_1,\dots, \sigma_{n-1}\, \begin{array}{|lr}
                                              \sigma_i\sigma_j=\sigma_j\sigma_i, & |i-j|>1 \\
                                                \sigma_i\sigma_{j}\sigma_i= \sigma_{j}\sigma_{i}\sigma_{j} , & |i-j|=1
                                             \end{array}
 \right\rangle .$$

\noindent The Coxeter group $W$ associated to $(A,\Sigma)$ can be obtained by adding the relations $\sigma_i^2=1$:
$$W=\langle \Sigma \,|\, \sigma_i^2=1 \, \forall i\in S ; \underbrace{\sigma_i\sigma_j\sigma_i\dots}_{m_{i,j} \text{ elements}}=\underbrace{\sigma_j\sigma_i\sigma_j\dots}_{m_{i,j} \text{ elements}} \forall i,j\in S,\, i\neq j,\, m_{i,j}\neq \infty \rangle.$$
If $W$ is finite, the corresponding Artin-Tits group is said to have spherical type. We will just consider Artin-Tits groups of spherical type, assuming that a spherical type Artin-Tits system is fixed. If $A$ cannot be decomposed as direct product of non-trivial Artin-Tits groups, we say that $A$ is \emph{irreducible}. Irreducible Artin-Tits groups of spherical type are completely classified \citep{Coxeter}.

\medskip
Let  $A$ be an Artin-Tits group of spherical type. A \textbf{standard parabolic subgroup}, $A_X$, is the subgroup generated by some $X\subseteq \Sigma$. A subgroup $P$ is called \textbf{parabolic} if it is conjugated to a standard parabolic subgroup, that is, $P=\alpha A_Y \alpha^{-1}$ for some standard parabolic subgroup $A_Y$ and some $\alpha\in A$. Notice that we may have $P=\alpha A_Y \alpha^{-1}=\beta A_Z \beta^{-1}$ for distinct $Y, Z\subset \Sigma$ and distinct $\alpha,\beta\in A$. We will write $P=(Y,\alpha)$ to express that $A_Y$ and $\alpha$ are known data defining the parabolic subgroup $P$. The conjugation of $P$ by $\alpha$ will also be noted by $P^{\alpha}$.  
\medskip

There is a natural way to associate a parabolic subgroup of $B_n$ to a curve system. Suppose that $A=B_n$ and let $A_X$ be the standard parabolic subgroup generated by $\{\sigma_i,\sigma_{i+1},\ldots , \sigma_j\}\subseteq \{\sigma_1,\ldots , \sigma_{n-1}\}$. Let $\mathcal{C}$ be the isotopy class of the circle enclosing the punctures $i,\ldots , j+1$ in $D_n$. Then $A_X$ fixes $\mathcal{C}$ and we will say that that $A_X$ is the parabolic subgroup associated to  $\mathcal{C}$. If $\mathcal{C}'=\mathcal{C}^\alpha$ for some $\alpha\in B_n$, then $(A_X )^\alpha:=\alpha^{-1}A_X\alpha$ is the parabolic subgroup associated to $\mathcal{C}'$. The parabolic subgroup associated to a system of non-nested curves is the direct sum of the subgroups associated to each curve. Notice that this is a well defined subgroup of $B_n$, as the involved subgroups commute. In this way, parabolic subgroups play a similar role, in Artin-Tits groups, to the one played by systems of curves in $B_n$.
\medskip

Our second purpose in this paper is to give a fast and simple algorithm to compute the minimal positive element that conjugates a given parabolic subgroup to a standard parabolic subgroup. The central Garside element of a standard parabolic subgroup $A_X$  will be denoted by $c_X$ and is to be defined in the next section. Having a generic parabolic subgroup, $P=(X, \alpha)$, the central Garside element will be denoted by $c_{P}$. We also define the minimal standardizer of the parabolic subgroup $P=(X, \alpha)$ to be the minimal positive element that conjugates $P$ to a standard parabolic subgroup. The existence and uniqueness of this element will be shown in this paper. Keep in mind that the $pn$-normal form of an element is a particular decomposition of the form $ab^{-1}$, where $a$ and $b$ are positive and have no common suffix. The main result of this paper is the following:

\begin{theorem} \label{minstand} Let $P=(X, \alpha)$ be a parabolic subgroup. If $c_{P}= ab^{-1}$ is in $pn$-normal form, then $b$ is the minimal standardizer of $P$.
\end{theorem}

Thus, the algorithm will take a parabolic subgroup $P=(X, \alpha)$ and will just compute the normal form of its central Garside element $c_{P}$, obtaining immediately the minimal standardizer of $P$. 

The paper will be structured in the following way: In section 2 some results and concepts about Garside theory will be recalled. In sections 3 and 4 the algorithm for braids will be explained. In section 5 the algorithm for Artin-Tits groups will be described and, finally, in section 6 we will bound the complexity of both procedures.

\section{Preliminaries about Garside theory}

Let us briefly recall some concepts from Garside theory. A group $G$ is called a \textbf{Garside group} with Garside structure $(G,\mathcal{P},\Delta)$ if it admits a submonoid $\mathcal{P}$ of positive elements such that $\mathcal{P}\cap \mathcal{P}^{-1}=\{1\}$ and a special element $\Delta \in \mathcal{P}$, called Garside element, with the following properties:

\begin{itemize}
 \item There is a partial order in $G$, $\po$,  defined by $a \po b \Leftrightarrow a^{-1}b \in \mathcal{P}$ such that for all $a,b\in G$ it exists a unique gcd $a \wedge b$ and a unique lcm $a \vee b$ with respect to $\po$.  This order is called prefix order and it is invariant by left-multiplication.

\item The set of simple elements $[1,\Delta]=\{a\in G\,|\, 1\po a \po \Delta  \}$ generates G.

\item $\Delta^{-1}\mathcal{P} \Delta = \mathcal{P}$. 

\item $\mathcal{P}$ is atomic: If we define the set of atoms  as the set of elements $a\in \mathcal{P}$ such that there is no non-trivial elements $b,c\in\mathcal{P}$ such that $a=bc$, then for every $x\in\mathcal{P}$ there is an upper bound on the number of atoms in the decomposition  $x=a_1a_2\cdots a_n$, where each $a_i$ is an atom.  

\end{itemize}

The conjugate by $\Delta$ of an element $x$ will be denoted $\tau(x)=x^\Delta =\Delta^{-1}x\Delta.$

In a Garside group, the monoid $\mathcal{P}$ also induces a partial order invariant under right-multiplication, the suffix order $\so$. This order is defined by $a\so b \Leftrightarrow ab^{-1}\in \mathcal{P}$, and for all $a,b\in G$ there exists an unique gcd $(a\wedger b)$ and an unique lcm $(a\veer b)$ with respect to $\so$. We say that a Garside group has \textbf{finite type} if $[1,\Delta]$ is finite. It is well known that Artin-Tits groups of spherical type admit a Garside structure of finite type \citep{Brieskorn1972,Dehornoy1999}. Moreover:

\begin{proposition}[\citealp{Vanderlek1983}]
A parabolic subgroup $A_X$ of an Artin-Tits group of spherical type is an Artin-Tits group of spherical type whose Artin-Tits system is $(A_X,X)$. 
\end{proposition}

\begin{proposition}[\citealp{Brieskorn1972}]
Let $(A_\Sigma,\Sigma)$ be an Artin-Tits system where $A_\Sigma$ is of spherical type. Then, a Garside element for $A_\Sigma$ is:
$$\Delta_\Sigma= \bigvee_{\sigma_i\in \Sigma}\sigma_i=\bigveer{\sigma_i\in \Sigma}\sigma_i,$$
and the submonoid of positive elements is the monoid generated by $\Sigma$.
Moreover, if $A_\Sigma$ is irreducible, then $(\Delta_\Sigma)^e$ generates the center of $A_\Sigma$, for some $e\in\{1,2\}$.
\end{proposition}

\begin{definition}
We define the right complement of a simple element $a$ as $\partial(a)=a^{-1}\Delta$ and the left complement as $\partial^{-1}(a)=\Delta a^{-1}$.
\end{definition}

\begin{remark} Observe that $\partial^2 =\tau$ and that, if $a$ is simple, then $\partial (a)$ is also simple, i.e., $1~\po~\partial(a)~\po~\Delta$. Both claims follow from $\partial(a)\tau(a)=\partial(a)\Delta^{-1}a\Delta = \Delta$ since $\partial(a)$ and $\tau(a)$ are positive.

\end{remark}

\begin{definition}
Given two simple elements $a,b$, the product $a\cdot b$ is said to be in left (resp. right) normal form if $ab\wedge \Delta =a$ (resp. $ab\wedger \Delta =b$). This is equivalent to say that $\partial(a)\wedge b =1$ (resp. $a\wedger \partial^{-1}(b)=1$).

We say that $x = \Delta^k x_1\cdots x_r$ is in \textbf{left normal form} if $k\in \mathbb{Z}$, $x_i\notin \{1,\Delta\}$ is a simple element for $i=1,\ldots , r$, and $x_i x_{i+1}$ is in left normal form for $0<i < r$. 

Analogously, $x =  x_1\cdots x_r\Delta^k$ is in \textbf{right normal form} if $k\in \mathbb{Z}$, $x_i\notin \{1,\Delta\}$ is a simple element for $i=1,\ldots , r$, and $x_i x_{i+1}$ is in right normal form for $0<i < r$. 

It is well known that the normal form of an element is unique \citep[Corollary 7.5]{Dehornoy1999}. Moreover, the numbers $r$ and $k$ do not depend on the normal form (left or right). We define the infimum, the canonical length and the supremum of $x$ respectively as  $\inf(x)=k$, $\ell(x)=r$ and $\sup(x)=k+r$.

\end{definition}

Let $a$ and $b$ be two simple elements such that $a\cdot b$ is in left normal form. One can write its inverse as $b^{-1}a^{-1}=\Delta^{-2}\partial^{-3}(b)\partial^{-1}(a)$. This is in left normal form because $\partial^{-1}(b)\partial(a)$ is in normal form by definition and $\tau=\partial^2$ preserves $\po$. More generally (see \citealp{Elrifai1994}), if $x=\Delta^k x_1\cdots x_r$ is in left normal form, then the left normal form of $x^{-1}$ is

$$x^{-1}= \Delta^{-(k+r)}\partial^{-2(k+r-1)-1}(x_r)\partial^{-2(k+r-2)-1}(x_{r-1})\cdots \partial^{-2k-1}(x_1)$$
For a right normal form, $x= x_1\cdots x_r \Delta^k$,  the right normal form of $x^{-1}$ is:

\begin{equation*}
x^{-1}= \partial^{2k+1}(x_r)\partial^{2(k+1)+1}(x_{r-1}) \cdots \partial^{2(k+r-1)+1}(x_1)\Delta^{-(k+r)} 
\end{equation*}

\begin{definition}[{\citealp[Theorem 2.6]{Charney1995}}]
Let $a,b \in \mathcal{P}$, then $x=a^{-1}b$ is said to be in $np$-normal form if $a\wedge b=1$. Similarly, we say that $x= ab^{-1}$ is in $pn$-normal form if $a\wedger b= 1$. 
\end{definition}

\begin{definition}\label{defirigid}Let  $\Delta^k x_1\cdots x_r$ with $r>0$ be the left normal form of $x$. We define the initial and the final factor respectively as $\iota(x)=\tau^{-k}(x_1)$ and $\varphi(x)=x_r$. We will say that $x$ is \textbf{rigid} if $\varphi(x)\cdot \iota(x)$ is in left normal form or if $r=0$. 
\end{definition}

\begin{definition}[{\citealp{Elrifai1994, Gebhardt2010a}}] Let $\Delta^k x_1\cdots x_r$ with $r>0$ be the left normal form of $x$. 
The \textbf{cycling} of $x$ is defined as $$\boldsymbol{c}(x)=x^{\iota(x)}=\Delta^k x_2\cdots x_r \iota(x). $$
The \textbf{decycling} of $x$ is $\boldsymbol{d}(x)=x^{(\varphi(x)^{-1})}= \varphi(x)\Delta^k x_1 \ldots x_{r-1}$.
We also define the preferred prefix of $x$ as $$\mathfrak{p}(x)= \iota(x) \wedge \iota(x^{-1}).$$
The \textbf{cyclic sliding} of $x$ is defined as the conjugate of $x$ by its preferred prefix: $$\mathfrak{s}(x)=x^{\mathfrak{p}(x)}=\mathfrak{p}(x)^{-1}x\mathfrak{p}(x).$$ 

\end{definition}

Let $G$ be a Garside group and $x\in G$. 
Keep in mind that $\inf_s(x)$ and $\sup_s(x)$ denote respectively the maximal infimum  and the minimal supremum in the conjugacy class $x^G$. 

\begin{itemize}

\item The super summit set of $x$  is 
\begin{align*}
SSS(x)=\,& \{y\in x^G \, |\, \ell \text{ is minimal in } x^G \} \\ 
 =\, & \{y\in x^G\,|\, \inf(y)={\inf}_s(y) \text{ and } \sup(y)={\sup}_s(y)\} 
\end{align*} 

\item The ultra summit set of $x$  is 
$$
\begin{array}{l}
 USS(x) =\{ y\in SSS(x) \,|\,  \boldsymbol{c}^m(y)=y  \text{ for some } m\geq 1\}  
\end{array} $$

\item The set of sliding circuits of $x$ is $$SC(x)=\{y\in x^G\,|\, \mathfrak{s}^m(y)=y \text{ for some } m\geq 1\}$$

\end{itemize}

These sets are finite if the set of simple elements is finite and their computation is very useful to solve the conjugacy problem in Garside groups. They satisfy the following inclusions: 

$$SSS(x)\supseteq USS(x) \supseteq SC(x).$$

\subsection*{The braid group, $B_n$}

A braid with $n$ strands can be seen as a collection of $n$ disjoint paths in a cylinder, defined up to isotopy, joining $n$ points at the top with $n$ points at the bottom, running monotonically in the vertical direction. 

\medskip

\begin{figure}[ht]
  \centering
  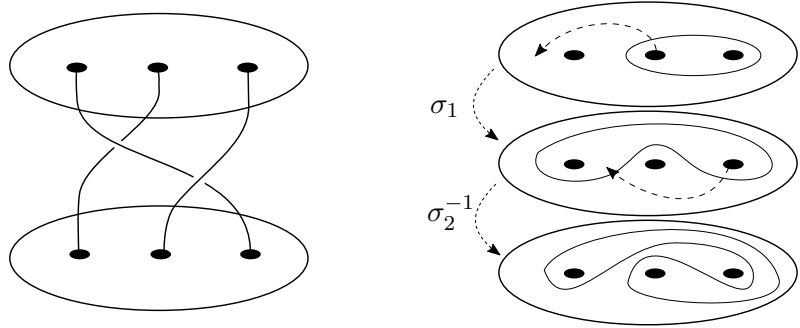
  \caption{The braid $\sigma_1\sigma_2^{-1}$ and how it acts on a curve in $D_3$.}
  \label{ct}
\end{figure}

Each generator $\sigma_i$ represents a crossing between the strands in positions $i$ and $i+1$ with a fixed orientation. The generator $\sigma_i^{-1}$ represents the same crossing with the opposite orientation. When considering a braid as a mapping class of $D_n$, these crossings are identified with the swap of two punctures in $D_n$ (See Figure~\ref{ct}).

\begin{remark}
%The braid group with $n$ strands, $B_n$, is a Garside group. Its generators are $\sigma_1,\dots,\sigma_{n-1}$, where each $\sigma_i$ represents a crossing between the strand $i$ and the strand $i+1$ and $\sigma_i^{-1}$ represent the same crossing with the opposite orientation. In this case 

After the above results, we see that the standard Garside structure of the braid group $B_n$ is $(B_n,B_n^+, \Delta_n)$ where
$$\Delta_n=\sigma_1\vee \dots \vee \sigma_{n-1}=(\sigma_1 \sigma_{2}\cdots \sigma_{n-1})(\sigma_1 \sigma_{2}\cdots \sigma_{n-2}) \cdots (\sigma_1 \sigma_2)\sigma_1$$ The simple elements in this case are also called \textbf{permutation braids} \citep{Elrifai1994}, because the set of simple braids is in one-to-one correspondence with the set of permutations of $n$ elements. 
Later we will use the following result: 
\end{remark}

\begin{lemma}[{\citealp[Lemma 2.4]{Elrifai1994}}]\label{strands}
If $s$ is a simple braid and its strands $j$ and $j+1$ cross, then $\sigma_j\po s$. 
\end{lemma}

\section{Detecting bending points} \label{sbending}

\begin{figure}[ht]
  \centering
  \includegraphics[width=0.3\linewidth]{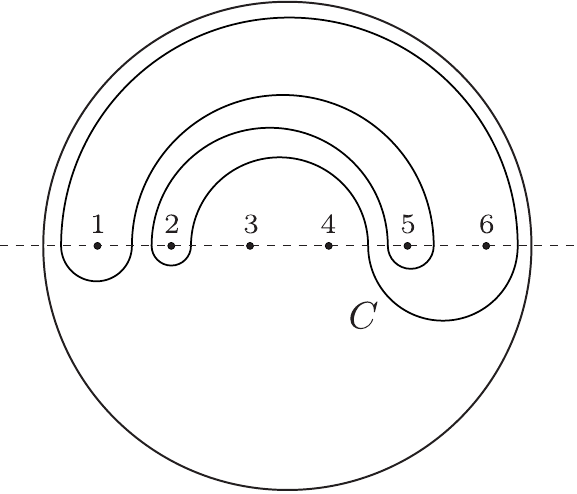}
  \caption{$W(C)=0\frown 6 \smile 4 \frown 2\smile 1\frown 4\smile 5 \frown 1\smile $.}
  \label{word}
\end{figure}

In order to describe a non-degenerate closed curve $C$ in $D_n$, we will use a notation introduced in \citep{Fenn1999}. Recall that $D_n$ has diameter $[0,n+1]$ and that the punctures of $D_n$ are placed at $1,2,\dots,n \in \mathbb{R}$. Choose a point on $C$ intersecting the real axis and choose an orientation for $C$. We will obtain a word $W(C)$ representing $C$, on the alphabet $\{\smile, \frown, 0,1,\dots, n\}$, by running along the curve, starting and finishing at the chosen point. We write down a symbol~$\smile$ for each arc on the lower half plane, a symbol $\frown$ for each arc on the upper half plane, and a number $m$ for each intersection of $C$ with the real segment $(m, m+1)$. An example is provided in Figure \ref{word}.  

\medskip

For an isotopy class of curves $\mathcal{C}$, $W(\mathcal{C})$ is the word associated to a reduced representative~$C^{red}$, i.e., a curve in $\mathcal{C}$ which has minimal intersection with the real axis. $C^{red}$ is unique up to isotopy of $D_n$ fixing the real diameter setwise \citep{Fenn1999}, and $W(\mathcal{C})$ is unique up to cyclic permutation and reversing. 

\begin{remark}
Notice that if a curve $C$ does not have minimal intersection with the real axis, then $W(C)$ contains a subword of the form $p\smile p\frown$ or $p\frown p\smile$. Hence, the curve can be isotoped by ``pushing'' this arc in order not to intersect the real axis. This is equivalent to remove the subword mentioned before from $W(C)$. In fact, we will obtain $W(\mathcal{C})$ by removing all subwords of this kind from $W(C)$. The process of removing $p\smile p\frown$ (resp. $p\frown p\smile$) from $W(C)$ is called \textbf{relaxation} of the arc $p\smile p$ (resp. $p\frown p$).
\end{remark}

\begin{figure}[ht]
  \centering
  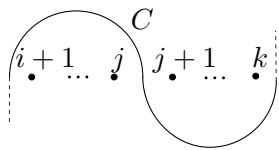
  \caption{A bending point at $j$ in a curve $C$.}
  \label{bending_point}
\end{figure}

\begin{definition}Let $C$ be a non-degenerate simple closed curve. We say that there is a \emph{bending point (resp. reversed bending point)} of $C$ at $j$ if we can find in $W(\mathcal{C})$, up to cyclic permutation and reversing, a subword of the form $i\frown j \smile k$ (resp. $i\smile j \frown k$) for some $0\leq i < j < k \leq n$ (Figure \ref{bending_point}). 

We say that a curve system has a bending point at $j$ if one of its curves has a bending point at $j$. 
\end{definition}

The algorithm we give in Section~\ref{salgo} takes a curve system $\mathcal{S}$ and `untangles' it in the shortest (positive) way. That is, it gives the shortest positive braid $\alpha$ such that $\mathcal{S}^\alpha$ is standard, i.e., the minimal element in $\St{S}$. Bending points are the key ingredient of the algorithm. We will show that if a curve system $\mathcal{S}$ has a bending point at $j$, then $\sigma_j$ is a prefix of the minimal element in $\St{S}$. This will allow to untangle $\mathcal{S}$ by looking for bending points and applying the corresponding $\sigma_j$ to the curve until no bending point is found. The aim of this section is to describe a suitable input for this algorithm and to show the following result.

\begin{proposition}\label{t1}
A curve system is standard if and only if its reduced representative has no bending points.
\end{proposition}

\subsection{Dynnikov coordinates}

\begin{figure}[ht]
  \bigskip
   \bigskip
  \centering
  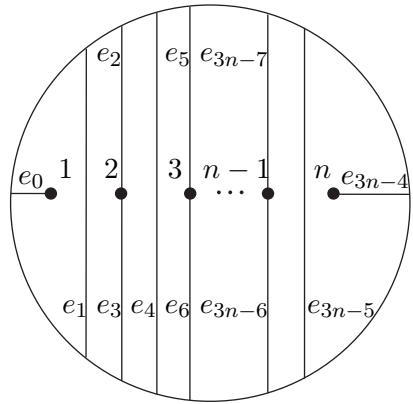
  \caption{Triangulation used to define Dynnikov coordinates.}
  \label{triangulation}
\end{figure}

\begin{figure}[h]
  \centering
  \includegraphics[width=0.5\linewidth]{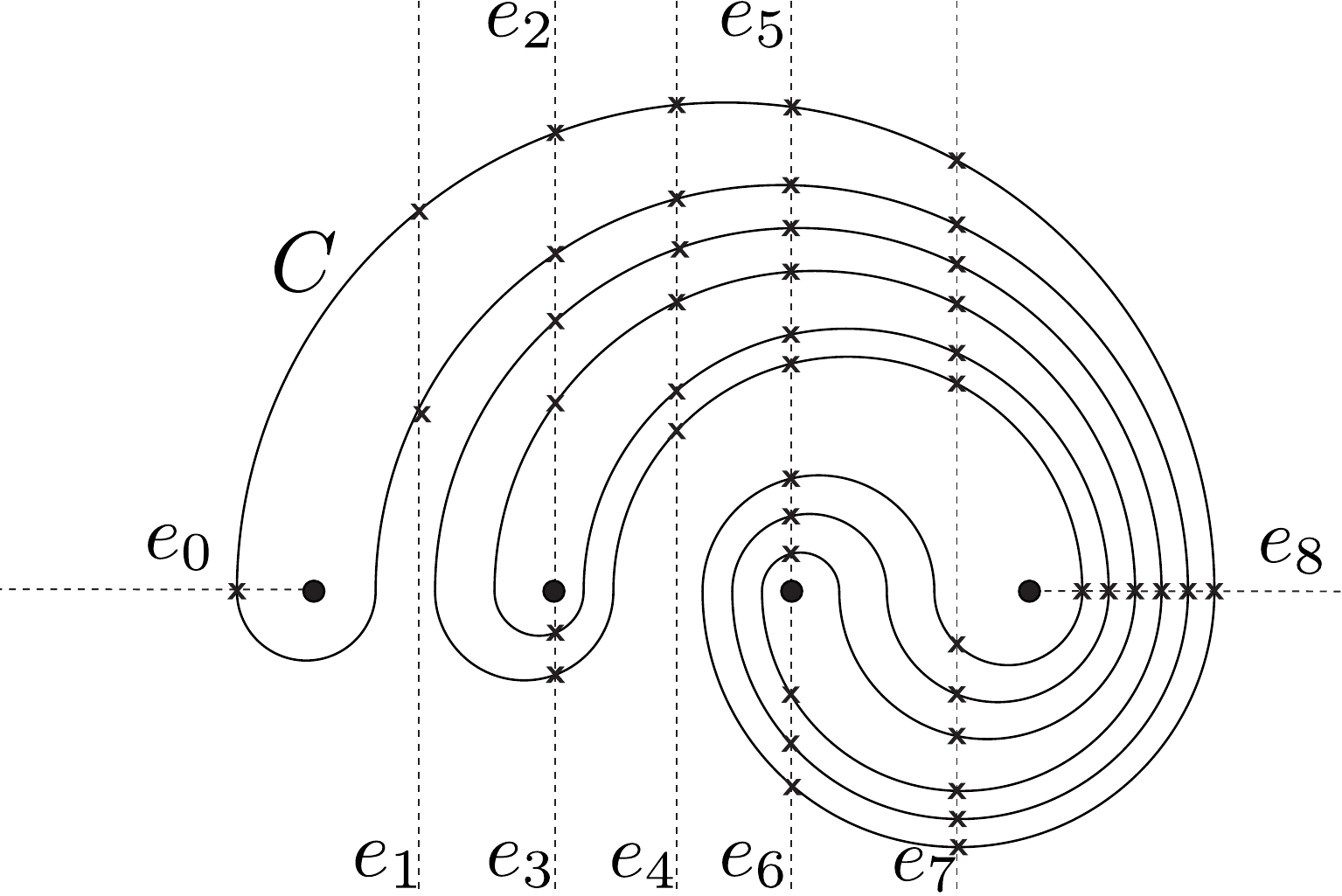}
  \caption{ The Dynnikov coordinates and reduced Dynnikov coordinates of $C$ are, respectively, $(x_0,\dots ,x_8)=(1,2,4,2,6,9,3,12,6)$;
  $(a_0,b_0, a_1,b_1,a_2,b_2, a_3, b_3)= (0,-1, 1,-2,3,-3,0, 6 ).$}
  \label{coord}
\end{figure}

We have just described a non-degenerate simple closed curve in $D_n$ by means of the word~$W(C)$. There is a different and usually much shorter way to determine a curve system $\mathcal{S}$ in~$D_n$: its Dynnikov coordinates \citep[Chapter~7]{Dehornoy2008}. The method to establish the coordinates of~$C$ is as follows. Take a triangulation of $D_n$ as in Figure \ref{triangulation} and let $x_i$ be the number of times the curve system $\mathcal{S}$ intersects the edge $e_i$. The Dynnikov coordinates of the curve system are given by the $t$-uple $(x_0,x_1,\dots , x_{3n-4})$. There exists a reduced version of these coordinates, namely $(a_0,b_0,\dots, a_{n-1},b_{n-1})$, where $$a_i = \dfrac{x_{3i-1} - x_{3i}}{2}, \qquad b_i = \dfrac{x_{3i-2} - x_{3i+1}}{2}, \qquad \forall i = 1, ... , {n-2}$$ and $a_0=a_{n-1}=0$, $b_0=-x_0$ and $b_{n-1}=x_{3n-4}$. See an example in Figure~\ref{coord}.

Furthermore, there are formulae determining how these coordinates change when applying $\sigma_j^{\pm 1}$, to the corresponding curve, for $0<j<n$.

\begin{proposition}[{\citealp[Proposition 8.5.4]{Dehornoy2002}}]\label{transformation} 
For $c=(a_0,b_0,\dots, a_{n-1},b_{n-1})$,  we have $$c^{\sigma_k^{-1}}=(a'_0,b'_0,\dots, a'_{n-1},b'_{n-1}),$$ with $a'_j=a_j, \, b'_j=b_j$ for $j\not\in \{k-1,k\}$,  and 

$$\begin{array}{rl}
a'_{k-1} & =a_{k-1}+(\delta^++b_{k-1})^+, \\ 
a'_{k} & =a_{k}-(\delta^+-b_{k})^+, \\ 
b'_{k-1} & =b_{k-1}-(-\delta')^++\delta^+, \\ 
b'_{k} & =b_{k}+(-\delta')^+-\delta^+,
\end{array}$$ 

where $\delta=a_{k}-a_{k-1}$, $\delta'=a'_{k}-a'_{k-1}$ and $x^+=\max(0,x)$. 

We also have $$c^{\sigma_k}=c^{\lambda\sigma_k^{-1}\lambda}$$ with $(a_1,b_1,\dots, a_{n-1},b_{n-1})^{\lambda}=(-a_1,b_1,\dots, -a_{n-1}, b_{n-1})$.

\end{proposition}

\begin{remark}
Notice that the use of $\sigma_k^{-1}$ in the first equation above is due to the orientation of the strands crossings that we are taking for our braids (see Figure~\ref{ct}), which is the opposite of the orientation used in \citep{Dehornoy2002}.
\end{remark}

\begin{figure}[ht]
\centering
\begin{subfigure}{.48\linewidth}
  \centering
  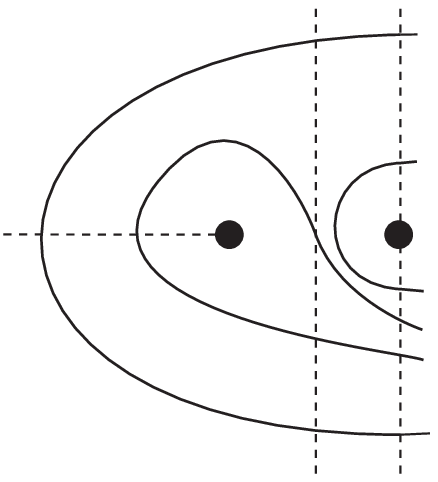
  \smallskip
  \caption{A bending point at 1.}
  \label{first_bending}
\end{subfigure}%
\begin{subfigure}{0.48\linewidth}
  \centering
  \bigskip
   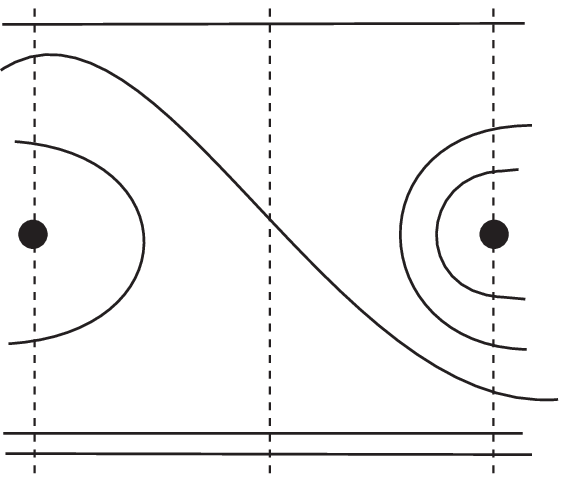
   \medskip
  \caption{A bending point and a right hairpin at $i$.}
  \label{middle_bending}
\end{subfigure}%
\caption{Detecting bending points with Dynnikov coordinates.}
\label{bend_coordinates}
\end{figure}

Let us see how to detect a bending point of a curve system $\mathcal{S}$ with these coordinates. Notice that there cannot be a bending point at $0$ or at $n$. It is easy to check that there is a bending point at $1$ if and only if $x_2<x_3$ (Figure \ref{first_bending}). Actually, if $R$ is the number of subwords of type  $0\frown 1 \smile k$ for some $1 < k \leq n$, then $x_3=x_2+2R$. Symmetrically, there is a bending point at $n-1$ if and only if $x_{3n-6}<x_{3n-7}$.

\bigskip

A bending point at $i$, for $1<i<n-1$, is detected by comparing the coordinates $a_{i-1}$ and~$a_i$ (Figure~\ref{middle_bending}). Notice that arcs not intersecting~$e_{3i-2}$ affect neither $a_{i-1}$ nor~$a_i$, and arcs not intersecting the real line do not affect the difference $a_{i-1}-a_{i}$. Hence, there is a bending point of~$\mathcal{S}$ at~$i$ if and only if $a_{i-1}-a_{i}>0$. Using a similar argument we can prove that there is a reversed bending point of~$\mathcal{S}$ at~$i$ if and only if $a_{i-1}-a_{i}<0$. Moreover, each bending point (resp. reversed bending point) at $i$ increases (resp. decreases) by 1 the difference $a_{i-1}-a_{i}$. We have just shown the following result:

\begin{lemma}[Bending point with Dynnikov coordinates] \label{bc}
Let $\mathcal{S}$ be a curve system on $D_n$ with reduced Dynnikov coordinates $(a_0,b_0,\dots, a_{n-1},b_{n-1})$.  For $j=1,\dots, n-1$ there are exactly $R$ \emph{bending points (resp. reversed bending points)} of $\mathcal{S}$ at~$j$ if and only if $a_{j-1}- a_{j}=R$ (resp. $a_{j-1}- a_{j}=-R$) .
\end{lemma}

\begin{lemma}\label{simetria}
Let $\mathcal{S}$ be a curve system as above. Then $\mathcal{S}$ is symmetric with respect to the real axis if and only if $a_i=0$, for $0<i<n$.
\end{lemma}

\emph{Proof.} Just notice that a symmetry with respect to the real axis does not affect $b$-coordinates and changes the sign of every $a_i$, for $0<i<n$. 
\hfill $\blacksquare$

\begin{lemma}\label{simetria2}
A curve system is standard if and only if it is symmetric with respect to the real axis. 
\end{lemma}

\emph{Proof}. 
For every $m=0,\ldots,n$, we can order the finite number of elements in $\mathcal S \cap (m,m+1)$ from left to right, as real numbers. Given an arc $a\frown b$ in $W(\mathcal S)$, suppose that it joins the $i$-th element in $\mathcal S \cap (a,a+1)$ with the $j$-th element in $\mathcal S\cap (b,b+1)$. The symmetry with respect to the real axis preserves the order of the intersections with the real line, hence the image $a\smile b$ of the above upper arc will also join the $i$-th element in $\mathcal S\cap (a,a+1)$ with the $j$-th element in $\mathcal S \cap (b,b+1)$. This implies that both arcs $a\frown b$ and $a \smile b$ form a single standard curve $a\frown b\smile$. As this can be done for every upper arc in $\mathcal S$, it follows that $\mathcal S$ is standard.
\hfill $\blacksquare$

\medskip

\emph{Proof of Proposition~\ref{t1}}. If the curve system is standard, then it clearly has no bending points. Conversely, if it has no bending points, by Lemma~\ref{bc} the sequence $a_0,\dots,a_{n-1}$ is non-decreasing, starting and ending at 0, so it is constant. By Lemmas~\ref{simetria} and \ref{simetria2}, the curve system is standard. \hfill $\blacksquare$

\section{Standardizing a curve system}\label{salgo}

\begin{figure}[ht]
  \bigskip
   \bigskip
  \centering
  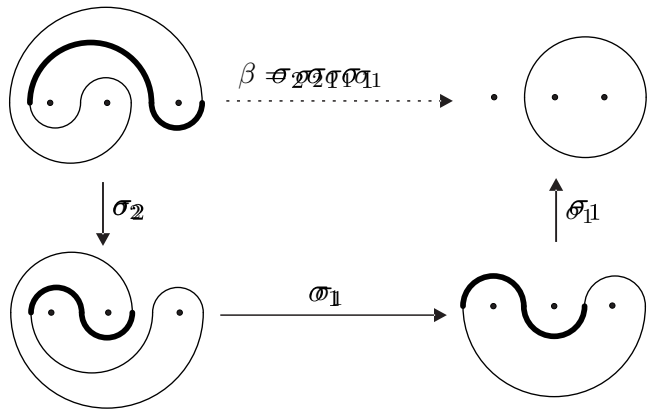
  \caption{A simple example of how to find the minimal standardizer of a curve.}
  \label{example}
\end{figure}

We will now describe an algorithm which takes a curve system $\mathcal{S}$, given in reduced Dynnikov coordinates, and finds the minimal element in $\St{S}$. The algorithm will do the following: Start with $\beta=1$. Check whether the curve has a bending point at $j$. If so, multiply $\beta$ by $\sigma_j$ and restart the process with $\mathcal{S}^{\sigma_j}$. A simple example is provided in Figure \ref{example}. The formal way is described in Algorithm \ref{algo}.

\begin{algorithm}
\SetKwInOut{Input}{Input}\SetKwInOut{Output}{Output}
\Input{The reduced coordinates $(a_0,b_0,\dots, a_{n-1},b_{n-1})$ of a curve system $\mathcal{S}$ on $D_n$.}
\Output{The $\po$-minimal element of $\St{S}$.}
\BlankLine
$c=(a_0,b_0,\dots, a_{n-1},b_{n-1})$;

$\beta= 1$;

$j=1$;

\While{$j<n$}{
    \If{$a_j<a_{j-1}$}{ 
      $c=c^{\sigma_j}$; (use Proposition~\ref{transformation})

      $\beta= \beta \cdot \sigma_j$;

      $j=1$;
      }
   \Else{$j=j+1$;}
}
 \Return{$\beta$;}
\caption{Rounding algorithm \label{algo}}
\end{algorithm}

\bigskip
The minimality of the output is guaranteed by the following theorem, which shows that $\sigma_j$ is a prefix of the minimal standardizer in $\St{S}$, provided $\mathcal{S}$ has a bending point at $j$.

\begin{theorem}\label{tprefix}
Let $\mathcal{S}$ be a curve system with a bending point at $j$.  Then $\sigma_j$ is a prefix of~$\alpha$, for every positive braid $\alpha$ such that $\mathcal{S}^\alpha$ is standard. 
\end{theorem} 

To prove the theorem we will need a result from \citep{Calvez2012}.

\begin{definition}
We will say that a simple braid $s$ is \emph{compatible} with a bending point at~$j$ if the strands $j$ and $j+1$  of $s$ do not cross in $s$. That is, if $\sigma_j\not\po s$.
\end{definition}

\begin{lemma}[{\citealp[Lemma 8]{Calvez2012}}]\label{calvez2}  Let $s_1$ and $s_2$ be two simple braids such that~$s_1 s_2$ is in left normal form. Let~$C$ be a curve with a bending point at~$j$ compatible with~$s_1$. Then, there exists some bending point of $C^{s_1}$ compatible with~$s_2$. 
\end{lemma}

\begin{remark}
The previous lemma holds also for a curve system, with the same proof.
\end{remark}

\emph{Proof of Theorem \ref{tprefix}.} Suppose that $\alpha \in B_n^{+}$ is such that $\mathcal{S}^\alpha$ is standard and $\sigma_j \not\preccurlyeq \alpha$.  Let~$s_1\cdots s_r$ be the left normal form of~$\alpha$. Notice that there is no $\Delta^p$ in the normal form. Otherwise, $\sigma_j$ would be a prefix of~$\alpha$, because it is a prefix of~$\Delta$. By Lemma~\ref{strands}, the strands~$j$ and $j+1$ of $s_1$ do not cross because $\sigma_j$ is not a prefix of~$\alpha$. Thus, $s_1$ is compatible with a bending point at~$j$ and by Lemma~\ref{calvez2}, $\mathcal{S}^{s_1}$~has a bending point compatible with~$s_2$. By induction, $\mathcal{S}^{s_1\cdots s_m }$~has a bending point compatible with~$s_{m+1}$, for $m=2,\dots,r$, where $s_{r+1}$ is the trivial braid (notice that Lemma~\ref{calvez2} is also valid if~$s_2$ is the trivial braid). Hence, $\mathcal{S}^{\alpha}$~has a bending point, i.e., it is not standard. A contradiction. \hfill $\blacksquare$

\medskip
In Algorithm \ref{algo} we can find the detailed procedure to compute the minimal element in $\St{S}$. Notice that, at every step, either the resulting curve has a bending point, providing a new letter of the minimal element in $\St{S}$, or it is standard and we are done. The process stops as $B_n^+$ has homogeneous relations, so all positive representatives of the minimal element have the same length, which is precisely the number of bending points found during the process. 

Notice that Theorem~\ref{tprefix} guarantees that the output of Algorithm~\ref{algo} is a prefix of every standardizer of $S$. This provides an alternative proof of the existence and uniqueness of a minimal element in $\St{S}$.

\section{Standardizing a parabolic subgroup}

Now we will give an algorithm to find the minimal standardizer of a parabolic subgroup $P=~(X, \alpha)$ of an Artin-Tits group $A$ of spherical type. The existence and uniqueness of this element will be shown by construction. 

%We will see that an element conjugate an parabolic subgroup, $P_1$, into another one, $P_2$, if and only if it conjugates the generator of the center of $P_1$ into the generator of the center of $P_2$. Then, we will show that to obtain the minimal standardizer of $P$ it suffices to compute the $pn$-normal form of the generator of its center.

\begin{definition}
Let $A_X$ be an Artin-Tits group of spherical type. We define its central Garside element as $c_X= (\Delta_X)^e$, where $e$ is the minimal positive integer such that $(\Delta_X)^e\in Z(A_X)$. We also define $c_{X,\alpha}:= \alpha c_X\alpha^{-1}$.
\end{definition}

\begin{proposition}[{\citealp[Proposition 2.1]{Godelle2003}}]\label{godequal}
Let $X,Y\subseteq \Sigma$ and $g\in A$. The following conditions are equivalent, 

\begin{enumerate}
\item $g^{-1}A_X g \subseteq A_Y$;
\item $g^{-1}c_X g \in A_Y $;
\item $g = xy$ where $y \in A_Y$ and $x$ conjugates $X$ to a subset of $Y$.
\end{enumerate}

\end{proposition}

The above proposition is a generalization of \citep[Theorem 5.2]{Paris1997} and implies, as we will see, that conjugating standard parabolic subgroups is equivalent to conjugating their central Garside elements. This will lead us to the definition of the central Garside element for a non-standard parabolic subgroup as given in Proposition~\ref{def_central_Garside}. In order to prove the following results, we need to define an object that generalizes to Artin-Tits groups of spherical type some operations used in braid theory:

%Artin-Tits groups of spherical type can be represented by Coxeter graphs. Recall that such a group, $A$, is defined by a symmetric matrix $M=(m_{i,j})_{i,j\in S}$  and the set of generators $\Sigma~=~\{\sigma_i\,|\, i\in S\}$, where $S$ is a finite set. The \emph{Coxeter graph} associated to $A$ is denoted $\Gamma_A$. The set of vertices of $\Gamma_A$ is $\Sigma$, and there is an edge joining two vertices $s,t\in \Sigma$ if $m_{s,t}\geq 3$. The edge will be labelled with $m_{s,t}$ if $m_{s,t}\geq 4$. 
%We say that the group $A$ is indecomposable if $\Gamma_A$ is connected and decomposable otherwise. If $A$ is decomposable, then there exists a non trivial partition $\Sigma=X_1 \cup\cdots \cup X_k$ such that $A$ is isomorphic to $A_{X_1}\times\cdots\times A_{X_k}$, where each $A_{X_j}$ is indecomposable. Each $A_{X_j}$ is called an \emph{indecomposable component} of $A$. 

\begin{definition}
Let $X\subset \Sigma$, $t\in \Sigma$.  We define 
$$r_{X,t}=\Delta_{X\cup \{t\}}\Delta^{-1}_{X}.$$
\end{definition}

\begin{remark}
In the case $t\notin X$, this definition is equivalent to the definition of positive elementary ribbon \citep[Definition 0.4]{Godelle2003}.  Notice that if $t\in X$, $r_{x,t}=1$.
\end{remark}

\begin{proposition}
There is a unique $Y\subset X\cup \{t\}$ such that $r_{X,t}X =Yr_{X,t}$.
\end{proposition}

\emph{Proof.}
Given $Z\subset \Sigma$, conjugation by $\Delta_Z$ permutes the elements of $Z$. Let us denote by $Y$ the image of $X$ under the permutation of $X\cup \{t\}$ induced by the conjugation by $\Delta_{X\cup \{t\}}$. Then
$$r_{X,t}X r_{X,t}^{-1} = \Delta_{X\cup \{t\}} \Delta_X^{-1} X \Delta_X \Delta_{X\cup \{t\}}^{-1} = \Delta_{X\cup \{t\}} X \Delta_{X\cup \{t\}}^{-1} = Y.$$ 
\hfill $\blacksquare$
 
% If $X,Y\subset \Sigma $ we say that $g \in A^{+}_\Sigma$ is a positive $Y$-ribbon-$X$ if $g =g_n \cdots g_1$ where $g_i$ is a positive elementary
%$X_i$-ribbon-$X_{i-1}$, with $X_0 =X$ and $X_n =Y$. 

\medskip

Artin-Tits groups of spherical type can be represented by Coxeter graphs. Recall that such a group, $A$, is defined by a symmetric matrix $M=(m_{i,j})_{i,j\in S}$  and the set of generators $\Sigma~=~\{\sigma_i\,|\, i\in S\}$, where $S$ is a finite set. The \emph{Coxeter graph} associated to $A$ is denoted $\Gamma_A$. The set of vertices of $\Gamma_A$ is $\Sigma$, and there is an edge joining two vertices $s,t\in \Sigma$ if $m_{s,t}\geq 3$. The edge will be labelled with $m_{s,t}$ if $m_{s,t}\geq 4$. We say that the group $A$ is indecomposable if $\Gamma_A$ is connected and decomposable otherwise. If $A$ is decomposable, then there exists a non trivial partition $\Sigma=X_1 \cup\cdots \cup X_k$ such that $A$ is isomorphic to $A_{X_1}\times\cdots\times A_{X_k}$, where each $A_{X_j}$ is indecomposable (each $X_j$ is just the set of vertices of a connected component of $\Gamma_X$). Each $A_{X_j}$ is called an \emph{indecomposable component} of $A$. 

\begin{lemma}\label{componentes}
Let $X,Y\subset \Sigma$ and let $X=X_1\sqcup \cdots \sqcup X_n$ and $Y=Y_1\sqcup \cdots \sqcup Y_m$ be the partitions of $X$ and $Y$, respectively, inducing the indecomposable components of $A_X$ and $A_Y$. Then, for every $g\in A$, the following conditions are equivalent:
\begin{enumerate}
\item $g^{-1} A_X g = A_Y$.

\item $m=n$ and $g=xy$, where $y\in A_Y$ and the parts of $Y$ can be reordered so that we have $x^{-1} {X_i} x~=~{Y_i}$ for $i=1,\ldots,n$.

\item $m=n$ and $g=xy$, where $y\in A_Y$ and the parts of $Y$ can be reordered so that we have $x^{-1} A_{X_i} x~=~A_{Y_i}$ for $i=1,\ldots,n$.
\end{enumerate}
\end{lemma}

\emph{Proof.}
Suppose that $g^{-1}A_X g = A_Y$. By  Proposition~\ref{godequal}, we can decompose $g=xy$ where $y\in A_Y$ and~$x$ conjugates the set $X$ to a subset of the set $Y$. Since conjugation by~$y$ induces an automorphism of $A_Y$, it follows that $x$ conjugates $A_X$ isomorphically onto~$A_Y$, so it conjugates $X$ to the whole set $Y$. Then conjugation by~$x$ sends indecomposable components of $X$ onto indecomposable components of $Y$. Hence $m=n$ and $x^{-1} {X_i} x= {Y_i}$ for $i=1,\ldots,n$ (reordering the indecomposable components of~$Y$ in a suitable way), as we wanted to show. Thus, statement 1 implies statement 2.

Now the relations satisfied by the elements of $X$ are also satisfied (through conjugation by~$x$) by their images in~$Y$, and viceversa. Since the connected components of~$\Gamma_X$ (resp. $\Gamma_Y$) are determined by the commutation relations among the letters of~$X$ (resp. $Y$), it follows that conjugation by~$x$ sends indecomposable components of $X$ onto indecomposable components of $Y$.Therefore statement 2 implies 3.

Finally, statement 3 implies 1 because $A_X=A_{X_1}\times \cdots \times A_{X_n}$ and ${A_Y=A_{Y_1}\times\cdots  \times A_{Y_n}}$.
\hfill $\blacksquare$

\begin{lemma}\label{precambiar}
Let $X,Y\subseteq \Sigma$, $g\in A$. Then, 
$$g^{-1}A_X g =A_Y \Longleftrightarrow g^{-1}c_X g =c_Y. $$

\end{lemma}

\emph{Proof.} Suppose that $g^{-1}c_X g =c_Y$. Then, by  Proposition~\ref{godequal}, we have $g^{-1}A_X g \subseteq A_Y$ and also $g A_Y g^{-1} \subseteq A_X$. As conjugation by $g$ is an isomorphism of $A$, the last inclusion is equivalent to $A_Y \subseteq  g^{-1}A_X g $. Thus, $g^{-1}A_X g =A_Y$, as desired. 

Conversely, suppose that $g^{-1}A_X g =A_Y $. By using Lemma~\ref{componentes}, we can decompose $g=xy$ where $y\in A_Y$ and $x$ is such that $x^{-1}A_{X_i} x =A_{Y_i}$, where $A_{X_i}$ and $A_{Y_i}$ are indecomposable components of~$A_X$ and $A_Y$ for $i=1,\dots, n$. As the conjugation by~$x$ defines an isomorphism between $A_{X_i}$ and $A_{Y_i}$, we have that $x^{-1}Z(A_{X_i}) x =Z(A_{Y_i})$. Hence, we have $x^{-1}c_{X_i} x =\Delta_{Y_i}^k$ for some $k\in\mathbb{Z}$. Let $c_{X_i}=\Delta_{X_i}^{e_1}$ and $c_{Y_i}=\Delta_{Y_i}^{e_2}$. As $A_{X_i}$ and $A_{Y_i}$ are isomorphic, then $e_1=e_2$. Also notice that in an Artin-Tits group of spherical type the relations are homogeneous and so $k=e_1=e_2$, having $x^{-1}c_{X_i} x =c_{Y_i}$. Let $$e=\max\{e_i\,|\, c_{X_i}=\Delta_{X_i}^{e_i}\}= \max\{e_i\,|\, c_{Y_i}=\Delta_{Y_i}^{e_i}\},$$ and denote $d_{X_i}=\Delta_{X_i}^e$ and $d_{Y_i}=\Delta_{Y_i}^e$ for $i=1,\ldots,n$. Notice that $d_{X_i}$ is equal to either $c_{X_i}$ or $(c_{X_i})^2$, and the same happens for each $d_{Y_i}$, hence $x^{-1} d_{X_i} x = d_{Y_i}$ for $i=1,\ldots, n$. Then, as $c_X=\prod_{i=1}^n{d_{X_i}}$ and $c_Y=\prod_{i=1}^n{d_{Y_i}}$, it follows that $x^{-1} c_X x = c_Y$. Therefore, $g^{-1} c_X g = y^{-1}(x^{-1}c_X x)y = y^{-1}c_Y y = c_Y.$
\hfill $\blacksquare$

\begin{lemma}\label{cambiar}
Let $P=(X,  \alpha)$ be a parabolic subgroup and $A_Y$ be a standard parabolic subgroup of an Artin-Tits group $A$ of spherical type. Then we have
$$g^{-1}P g=A_Y \Longleftrightarrow g^{-1}c_{X,\alpha} g = c_Y.$$

\end{lemma}

\emph{Proof}.
If $P=(X, \alpha)$, it follows that $g^{-1}Pg=A_Y$ if and only if $g^{-1}\alpha A_X \alpha^{-1}g=A_Y$. By Lemma~\ref{precambiar}, this is equivalent to $g^{-1}\alpha c_X \alpha^{-1} g =c_Y$, i.e., $g^{-1}c_{X,\alpha} g = c_Y$. \hfill $\blacksquare$

\begin{proposition}\label{def_central_Garside}
Let $P=(X,\alpha)=(Y,\beta)$ be a parabolic subgroup of an Artin-Tits group of spherical type. Then $c_{X,\alpha}=c_{Y,\beta}$ and we can define $c_P:=c_{X,\alpha}$ to be the central Garside element of $P$.
\end{proposition}

\emph{Proof.}
Suppose that $g$ is a standardizer of $P$ such that  $g^{-1}P g=A_Z$. By using Lemma~\ref{cambiar}, we have that $c_Z=g^{-1}c_{X,\alpha} g = g^{-1}c_{Y,\beta} g $. Thus, $c_{X,\alpha}=c_{Y,\beta}$. \hfill $\blacksquare$

\bigskip
By Lemma~\ref{cambiar}, to find the minimal standardizer of a parabolic subgroup~$P=(X,\alpha)$, we just need to find the minimal positive element conjugating $c_{P}$ to some $c_Y$. Let $$C^{+}_{A_\Sigma}(c_{P})=\{s\in \mathcal{P}\,|\, s =u^{-1}c_{P}u,\, u\in A_\Sigma \}$$ be the set of positive elements conjugate to~$c_{P}$ (which coincides with the positive elements conjugate to~$c_X$).  Firstly, let us compute the minimal conjugator from~$c_{P}$ to this set. That is, the shortest positive element $u$ such that $u^{-1}c_{X,\alpha}u \in \mathcal{P}$.

\begin{proposition}[{\citealp[Proposition 4.8]{Franco2003}}]\label{franco}
For any $\alpha\in A_\Sigma$ conjugated to a positive element, it exists a unique $\po$-minimal positive element conjugating $\alpha$ to $C^{+}_{A_\Sigma}(\alpha)$.
\end{proposition}

\begin{proposition}\label{minconj}
If $x= ab^{-1}$ is in $pn$-normal form and $x$ is conjugated to a positive element, then $b$ is a prefix of every positive element conjugating~$x$ to $C^{+}_{A_\Sigma}(x)$.
\end{proposition}

\emph{Proof.} Suppose that $\rho$ is a positive element such that $\rho^{-1}x\rho$ is positive. Then $1\po \rho^{-1}x\rho$. Multiplying from the left by $x^{-1}\rho$ we obtain $x^{-1}\rho\po \rho$ and, since $\rho$ is positive, $x^{-1}\po x^{-1}\rho\po \rho$. Hence $x^{-1}\po \rho$ or, in other words $ba^{-1}\po \rho$. On the other hand, by definition of $pn$-normal form, we have $a\wedger b=1$, which is equivalent to $a^{-1}\vee b^{-1}=1$ \citep[Lemma 1.3]{Gebhardt2010}. Multiplying from the left by $b$, we obtain $ba^{-1}\vee 1=b$.

Finally, notice that $ba^{-1}\po\rho$ and also $1\po \rho$. Hence $b=ba^{-1}\vee 1 \po \rho$. Since $b$ is a prefix of $\rho$ for every positive $\rho$ conjugating $x$ to a positive element, the result follows.
\hfill  $\blacksquare$

\begin{lemma}\label{ribbonsufijo}
Let $A_X$ be a standard parabolic subgroup and $t\in \Sigma$. If $\alpha\Delta_X\so t$, then $\alpha \so r_{X,t}$.

\end{lemma}
\emph{Proof.}
Since the result is obvious for $t\in X$ ($r_{X,t}=1$), suppose $t\notin X$. Trivially,  $\alpha\Delta_X\so \Delta_X$. As $\alpha\Delta_X\so t$, we have that $\alpha\Delta_X \so~\Delta_X~\veer ~t$. By definition, $\Delta_X \veer t~=~\Delta_{X\cup\{t\}}=r_{X,t}\Delta_X$. Thus, $\alpha\Delta_X\so r_{X,t} \Delta_X$ and so $\alpha \so r_{X,t}$, because $\so$ is invariant under right-multiplication.

\hfill$\blacksquare$

\begin{theorem*}
Let $P=(X, \alpha)$ be a parabolic subgroup. If $c_{P}= ab^{-1}$ is in $pn$-normal form, then $b$ is the $\po$-minimal standardizer of $P$.

\end{theorem*}

\emph{Proof.} 
We know from Proposition~\ref{minconj} that $b$ is a prefix of any positive element conjugating $c_P$ to a positive element, which guarantees its $\po$- minimality. We also know from Lemma~\ref{cambiar} that any standardizer of $P$ must conjugate $c_{P}$ to a positive element, namely to the central Garside element of some standard parabolic subgroup. So we only have to prove that $b$ itself conjugates $c_P$ to the central Garside element of some standard parabolic subgroup. We assume $\alpha$ to be positive, because there is always some $k\in \mathbb{N}$ such that $\Delta^{2k} \alpha$ is positive and, as $\Delta^2$ lies in the center of $A$, $P=(X,\alpha)=(X,\Delta^{2k}\alpha)$.

The $pn$-normal form of $c_{P}=\alpha c_X \alpha^{-1}$ is obtained by cancelling the greatest common suffix of $\alpha c_X$  and $\alpha$. Suppose that $t\in \Sigma$ is such that $ \alpha\so t$ and $\alpha c_X \so t$. 

If $t\notin X$, then $r_{X,t}\neq 1$ and by Lemma~\ref{ribbonsufijo} we have that $\alpha \so r_{X,t} $, i.e., $\alpha=\alpha_1 r_{X,t}$ for some $\alpha_1\in A_\Sigma$. Hence, 
$$\alpha c_X \alpha^{-1}=\alpha_1 r_{X,t} c_X r_{X,t}^{-1}\alpha_1^{-1}=\alpha_1 c_{X_1} \alpha_1^{-1}$$ 
for some $X_1\subset \Sigma$. In this case, we reduce the length of the conjugator (by the length of $r_{X,t}$). If $t\in X$, $t$ commutes with $c_X$, which means that
$$\alpha c_X \alpha^{-1}=\alpha_1 t c_X t^{-1}\alpha_1^{-1}=\alpha_1 c_{X_1} \alpha_1^{-1},$$ where $\alpha_1$ is one letter shorter than $\alpha$ and $X_1=X$.

We can repeat the same procedure for $\alpha_i c_{X_i} \alpha_i^{-1}$, where $X_i\subset \Sigma$, $t_i\in \Sigma$ such that $\alpha_i\so t_{i} $ and $\alpha_i c_{X_i} \so t_{i}$. As the length of the conjugator decreases at each step, the procedure must stop, having as a result the $pn$-normal form of $c_{P}$, which will have the form:
$$c_{P}= (\alpha_k c_{X_k}) \alpha_k^{-1},\quad \text{for } k\in\mathbb{N},\, X_k\subset \Sigma.$$

Then, $\alpha_k=b$ clearly conjugates $c_{P}$ to $c_{X_k}$, which is the central Garside element of a standard parabolic subgroup, so $b$ is the $\po$-minimal standardizer of $P$.
\hfill $\blacksquare$

%
%
%\begin{lemma}\cite[Lemma 5.6]{Paris1997}
%Let $u\in A_\Sigma^+$, $X\subseteq\Sigma$ and $\Delta_X^e$ the generator of $Z(A_X)$. If $u^{-1}\Delta_X^e u $ is positive, then 
%$$\exists Y\subseteq \Sigma \text{ such that } u^{-1}\Delta^{e}_Xu=\Delta^{e}_Y.$$ 
%
%\end{lemma}
\medskip
We end this section with a result concerning the conjugacy classes of elements of the form $c_{P}$. As all the elements of the form $c_Z$, $Z\subseteq X$, are rigid (Definition~\ref{defirigid}), using the next theorem we can prove that the set of sliding circuits of $c_{P}$ is equal to its set of positive conjugates.

\begin{theorem}[{\citealp[Theorem 1]{Gebhardt2010a}}]\label{rigid}
Let $G$ be a Garside group of finite type. If $x\in G$ is conjugate to a rigid element, then $SC(x)$ is the set of rigid conjugates of $x$.
\end{theorem}

\begin{corollary}\label{deltaconjuntos} Let $P=(X, \alpha)$ be a parabolic subgroup of an Artin-Tits group of spherical type. Then
$$C^{+}_{A_\Sigma}(c_{P})=SSS(c_P)=USS(c_P)=SC(c_P)= \{c_Y \,|\, Y\in \Sigma, \,c_Y \text{ conjugate to } c_X\}.$$
\end{corollary}

\emph{Proof.} By Theorem~\ref{rigid}, it suffices to prove that $C^{+}_{A_\Sigma}(c_P)$ is composed only by rigid elements of the form $c_Z$. Let $P'=(X,\beta)$ and suppose that $c_{P'}\in C^{+}_{A_\Sigma}(c_P)$. Let $b$ be the minimal standardizer of $c_{P'} $. By Proposition~\ref{minconj} and Theorem~\ref{minstand}, $b$ is the minimal positive element conjugating $c_{P'}$ to $C^{+}_{A_\Sigma}(c_{P})$, which implies that $b=1$, so $P'$ is standard. Hence, all positive conjugates of $c_{P'}$ are equal to $c_{Y}$ for some $Y$, therefore they are rigid. 
 \hfill $\blacksquare$

\begin{corollary}
Let $P=(X, \alpha)$ be a parabolic subgroup of an Artin-Tits group of spherical type. Then the set of positive standardizers of $P$, $$\mathrm{St}(P)=\{\alpha \in A_\Sigma^{+},\,|\, c_P^\alpha=c_Y, \, \text{for some } Y\subseteq \Sigma\},$$
is a sublattice of $A_\Sigma^{+}$.
\end{corollary}

\emph{Proof.} Let $s_1$ and $s_2$ be two positive standardizers of $P$ and let $\alpha := s_1 \wedge s_2$ and ${\beta := s_1 \vee s_2}$. By Corollary~\ref{deltaconjuntos} and, for example, \citep[Proposition~7, Corollary~8]{Gebhardt2010a}, we have that $c_P^\alpha =c_Y$ and $c_P^\beta =c_Z$ for some $Y,Z\subseteq \Sigma$. Hence $\alpha,\beta \in \mathrm{St}(P)$, as we wanted to show.
\hfill $\blacksquare$

\section{Complexity}

In this section we will describe the computational complexity of the algorithms which compute minimal standardizers of curves and parabolic subgroups. Let us start with Algorithm~\ref{algo}, which computes the minimal standardizer of a curve system.

Notice that Algorithm~\ref{algo} takes at each step the leftmost bending point of the curve system. However, Theorem~\ref{tprefix} allows us to take any bending point. The complexity of Algorithm~\ref{algo} will depend on the length of the output, which is the number of steps of the algorithm. To bound this length, we will compute a positive braid which belongs to $\St{S}$. This will bound the length of the minimal standardizer of $\mathcal{S}$.

\medskip

\begin{figure}[ht]
  \centering
  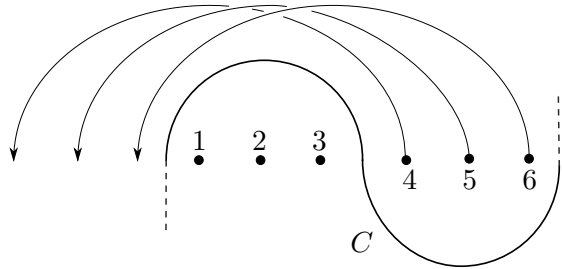
  \caption{Applying $s(0,3,6)$.}
  \label{eliminate}
\end{figure}

The usual way to describe the length (or the complexity) of a curve system consists in counting the number of intersections with the real axis, i.e., $\ell(\mathcal{S})= \#(\mathcal{S}\cap \mathbb{R})$. For integers $0\leq i < j < k\leq n$, we define the following braid (see Figure \ref{eliminate}):

$$s(i,j,k)=(\sigma_j\sigma_{j-1}\cdots \sigma_{i+1}) (\sigma_{j+1}\sigma_{j}\cdots \sigma_{i+2})\cdots(\sigma_{k-1}\sigma_{k-2}\cdots \sigma_{i+k-j})$$

\begin{lemma}\label{s} Applying $s=s(i,j,k)$ to a curve system $\mathcal{S}$, when $i\frown j \smile  k$ is a bending point, decreases the length of the curve system at least by two. 
\end{lemma}

\emph{Proof.}
We will describe the arcs of the curves of $\mathcal{S}$ in a new way, by associating a real number $c_p\in (0, n+1)$ to each of the intersections of $\mathcal{S}$ with the real axis, where $p$ is the position of the intersection with respect to the other intersections:  $c_1$ is the leftmost intersection and $c_{\ell(\mathcal{S})}$ is the rightmost one.  We will obtain a set of words  representing the curves of $\mathcal{S}$, on the alphabet $\{\smallsmile, \smallfrown, c_1,\dots, c_{\ell(\mathcal{S})}\}$, by running along the curve, starting and finishing at the same point. As before, we write down a symbol~$\smallsmile$ for each arc on the lower half plane, and a symbol $\smallfrown$ for each arc on the upper half plane. We also define the following function that sends this alphabet to the former one: $$L\colon \{\smallsmile, \smallfrown, c_1,\dots, c_{\ell(\mathcal{S})}\}\longrightarrow \{\smile, \frown, 0,\dots, n\}$$ $$L(\smallsmile)=\smile,\quad L(\smallfrown)=\frown,\quad  L(c_p)= \lfloor c_p \rfloor. $$ 
Take a disk $D$ such that the $\partial(D)$ intersects the real axis at two points, $x_2$ and $x_3$. Consider another point $x_1$ on the real axis such that $L(x_1)<L(x_2)$. Suppose that there are no arcs of $\mathcal{S}$ on the upper-half plane intersecting the arc $x_1\smallfrown x_2$ and there are no arcs of $\mathcal{S}$ on the lower-half plane intersecting the arc $x_2\smallsmile x_3$. We denote $I_1=(0, x_1)$, $I_2=(x_1,x_2)$, $I_3=(x_2,x_3)$ and $I_4=(x_3,n+1)$ and define $|I_t|$ as the number of punctures that lie in the interval $I_t$. 

\begin{figure}[ht]
\centering

\begin{subfigure}{.48\linewidth}
  \centering
  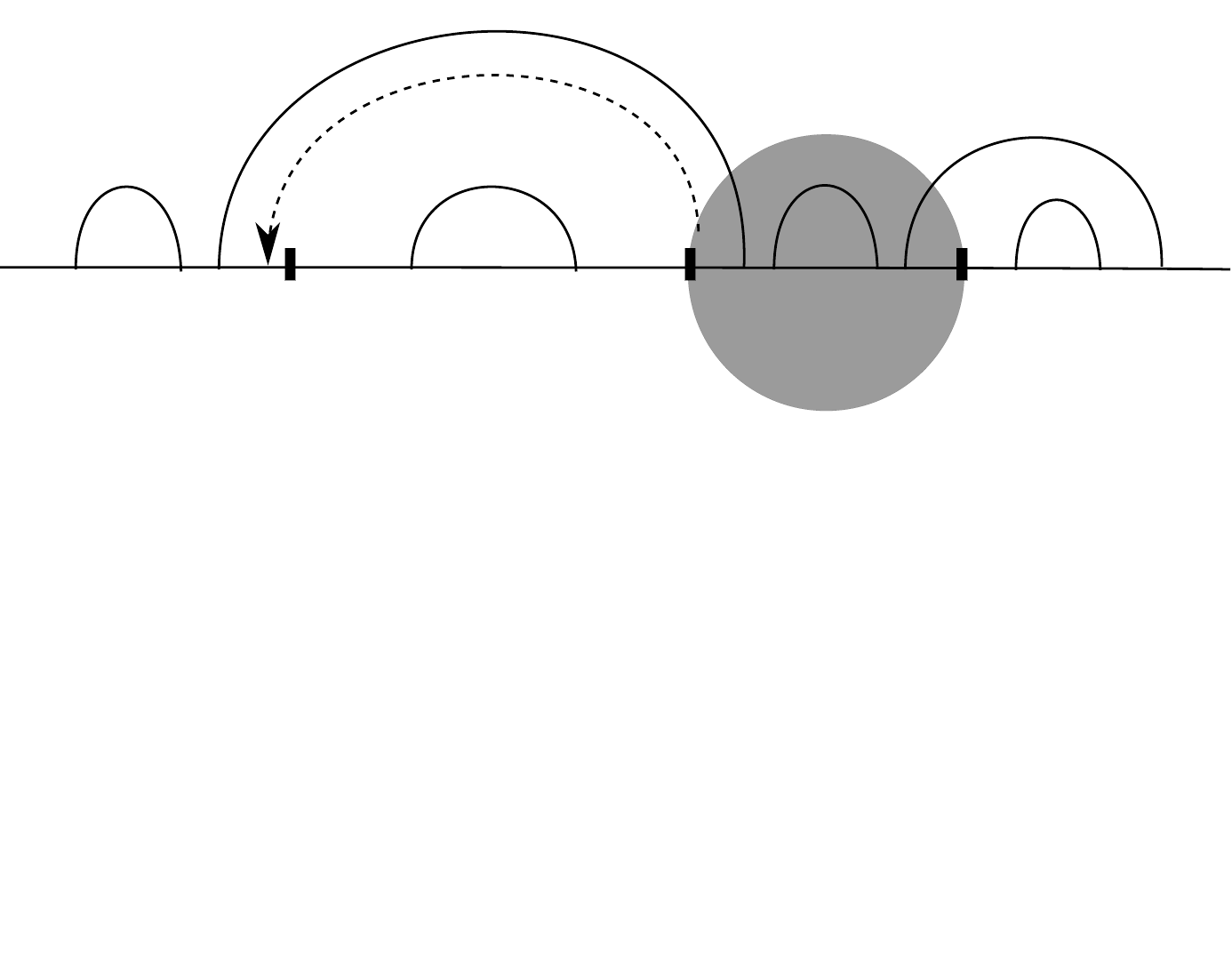
  \caption{$d$ acting on the arcs in the upper half plane}
  \label{arcossup}
\end{subfigure}%
\begin{subfigure}{0.48\linewidth}
  \centering
   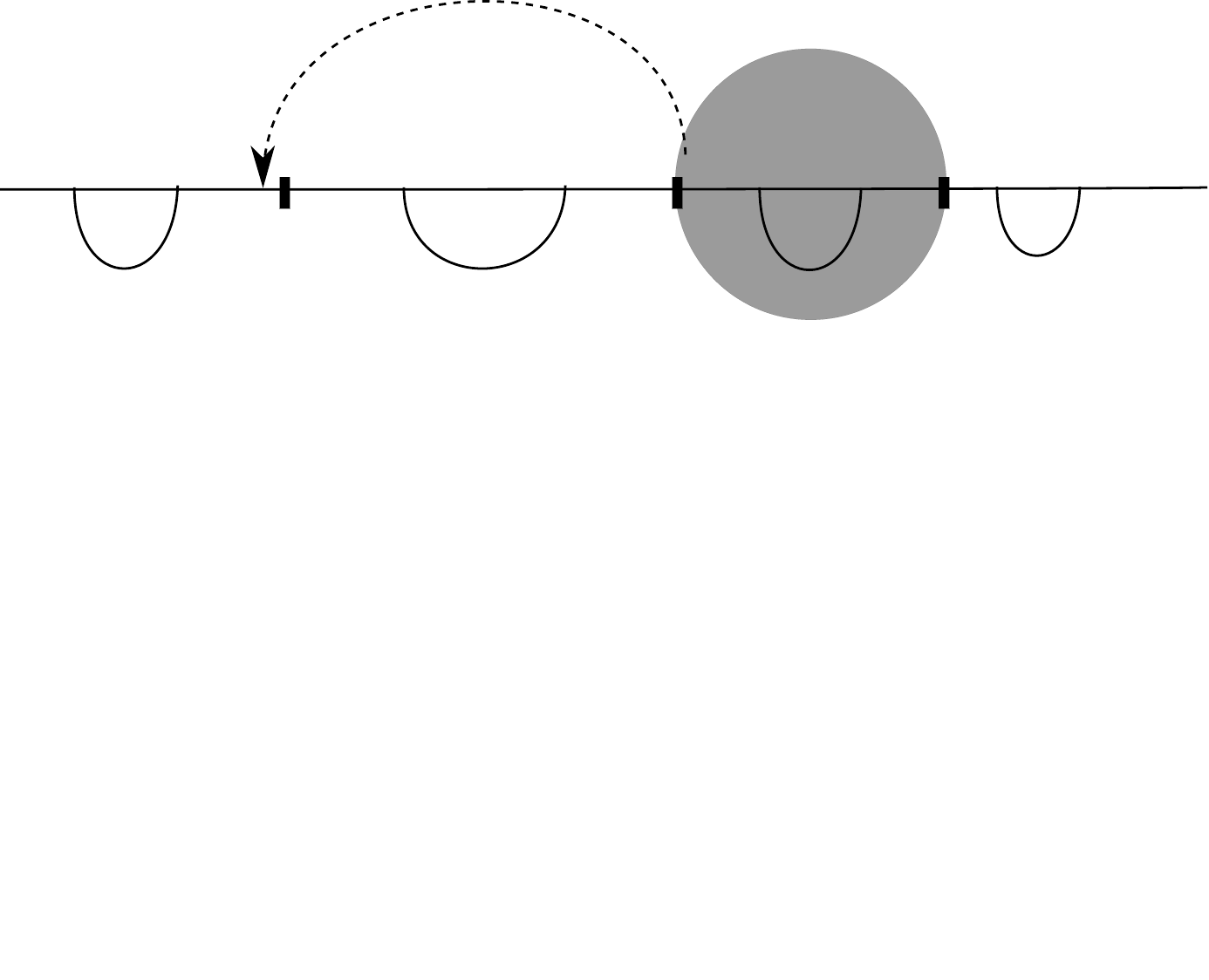
  \caption{$d$ acting on the arcs in the lower half plane}
  \label{arcosinf}
\end{subfigure}%
\caption{How the automorphism $d(x_1,x_2,x_3)$ acts on the arcs of $C$.}
\label{Arcos}
\end{figure}

We consider an automorphism of $D_n$, called $d=d(x_1,x_2,x_3)$, which is the final position of an isotopy that takes $D$ and moves it trough the upper half-plane to a disk of radius $\epsilon$ centered at $x_1$, which contains no point $c_p$ and no puncture, followed by an automorphism which fixes the real line as a set and takes the punctures back to the positions $1,\dots , n$. This corresponds to ``placing the interval $I_3$ between the intervals $I_1$ and $I_2$".  Firstly, we can see in Figure~\ref{Arcos} that the only modifications that  the arcs of $\mathcal{S}$ can suffer is the shifting of their endpoints. By hypothesis, there are no arcs in the upper half-plane joining $I_2$ with $I_j$ for $j\neq 2$, and there are no arcs in the lower half-plane joining $I_3$ with $I_j$ for $j\neq 3$. Any other possible arc is transformed by $d$ into a single arc, so every arc is transformed in this way. Algebraically, take an arc of $\mathcal{S}$, $c_{a}\smallfrown c_{b}$ (resp. $c_{a}\smallsmile c_{b}$ ), such that $L(c_{a})=\tilde{a}$ and $L(c_{b})=\tilde{b}$. Then, its image under $d$ is $c'_{a}\smallfrown c'_{b}$ (resp. $c'_{a}\smallsmile c'_{b}$ ) where

$$L(c'_p)=\left\{\begin{array}{ll}
\tilde{p} & \text{if } c_p\in I_1,I_4, \\ 
\tilde{p} + |I_3| & \text{if } c_p\in I_2, \\ 
\tilde{p} -|I_2| & \text{if } c_p\in I_3,
\end{array}   \right. \quad \text{ for } p=a,b.$$

After applying this automorphism, the curve could not be reduced, so relaxation of unnecessary arcs could be done, reducing the complexity of $\mathcal{S}$.

\medskip

\begin{figure}[ht]
  \centering
  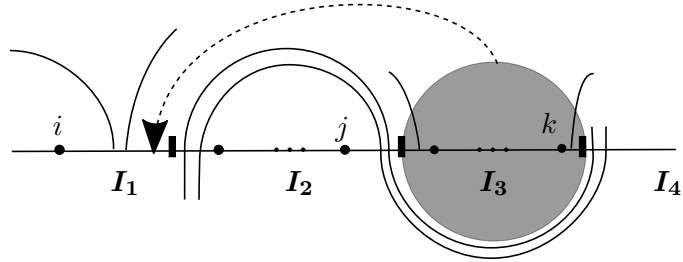
  \caption{Applying $s(i,j,k)$ to a curve is equivalent to permute their intersections with the real axis and then make the curve tight.}
  \label{homeo}
\end{figure}

Now, given a bending point $i\frown j\smile k$ of $\mathcal{S}$, consider the set $$B=\{c_p\smallfrown c_q \smallsmile c_r \,|\, L(c_p)<L(c_q)<L(c_r) \text{ and } L(c_q)=j\}$$ and choose the element of $B$ with greatest subindex $q$, which is also the one with lowest $p$ and $r$. Define $x_1,x_2$ and $x_3$ such that $x_1\in (c_{p-1},c_p)\cap (L(c_p),L(c_p)+1)$, $x_2\in (c_q,c_{q+1})\cap (j,j+1)$ and $x_3\in (c_{r-1},c_r)\cap (L(c_r),L(c_r)+1)$. Then, the braid $s(L(c_p), j, L(c_r))$ is represented by the automorphism $d(x_1,x_2,x_3)$ (see Figure~\ref{homeo}). Notice that the choice of the bending point from $B$ guarantees the non-existence of arcs of $C$ intersecting $x_1\smallfrown x_2$ or $x_2\smallsmile x_3$. After the swap of $I_2$ and $I_3$, the arc $c_q\smallsmile c_r$ will be transformed into $c'_q\smile c'_r$, where $L(c'_q)=L(c'_r)=L(c_r)$, and then relaxed, reducing the length of $\mathcal{S}$ at least by two. \hfill $\blacksquare$

 \medskip 
 
The automorphism $s=s(i,j,k)$ involves at most $(n-1)^2$ generators. Then, the output of our algorithm has at most $\frac{1}{2}\ell(\mathcal{S})(n-1)^2$ letters, because we have proven that $s$ reduces the length of the curve system at each step. Let us bound this number in terms of the input of the algorithm, i.e., in terms of reduced Dynnikov coordinates.

\begin{figure}[ht]
\centering
\begin{subfigure}{.48\linewidth}
  \centering
  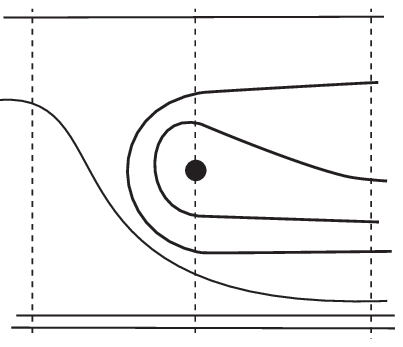
  \medskip
  \caption{Two left hairpins}
  \label{hairpin_left}
\end{subfigure}%
\begin{subfigure}{0.48\linewidth}
  \centering
   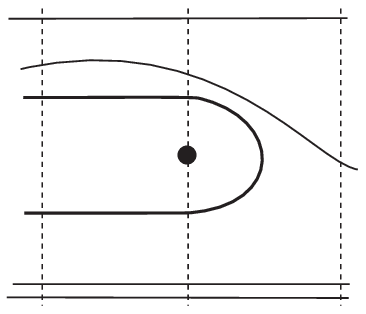
  \smallskip
  \caption{A right hairpin}
  \label{hairpin_right}
\end{subfigure}%
\caption{Detecting hairpins with Dynnikov coordinates.}
\label{bend_coordinates_hairpin}
\end{figure}

\begin{definition}
We say that there is a \emph{left hairpin} (resp. a \emph{right hairpin}) of $C$ at $j$ if we can find in $W(\mathcal{C})$, up to cyclic permutation and reversing, a subword of the form $i\frown j-1 \smile k$ (resp. $i\frown j \smile k$)  for some $i,k > j-1 $ (resp. $i,k<j$) (see Figure~\ref{bend_coordinates_hairpin}).
\end{definition}

\begin{proposition}\label{cotaminimalstand}
Let $\mathcal{S}$ be a curve system on $D_n$ represented by the reduced Dynnikov coordinates $(a_0,b_0,\dots, a_{n-1},b_{n-1})$. Then $\ell(\mathcal{S})\leq\sum_{i=0}^{n-1} (2|a_i| + |b_i|)$.
\end{proposition}

\emph{Proof.} Notice that each intersection of a curve $\mathcal{C}$ with the real axis corresponds to a subword of $W(\mathcal{C})$ of the form $i \frown j \smile k$ or $i\smile j \frown k$. If $i<j<k$ the subword corresponds to a bending point or a reversed bending point respectively. If $i,k>j$, there is a left hairpin at $j+1$. Similarly, if $i,k< j$, there is a right hairpin at $j$.

Recall that Lemma~\ref{bc} already establishes how to detect bending points with reduced Dynnikov coordinates. In fact, there are exactly $R$ bending points (including reversed ones) at $i$ if and only if $|a_{i-1}- a_{i}|=R$. We want to detect also hairpins in order to determine $\ell(\mathcal{S})$. Observe in Figure~\ref{bend_coordinates_hairpin} that the only types of arcs that can appear in the region between the lines $e_{3j-5}$ and $e_{3j-2}$ are left or right hairpins and arcs intersecting both $e_{3j-5}$ and $e_{3j-2}$. The arcs intersecting both $e_{3j-5}$ and $e_{3j-2}$ do not affect the difference $x_{3j-5}-x_{3j-2}$ whereas each left hairpin decreases it by 2 and each right hairpin increases it by 2. Notice that in the mentioned region there cannot be left and right hairpins at the same time. Then, there are exactly $R$ left (resp. right) hairpins at $j$ if and only if $b_{j-1}=-R$ (resp. $b_{j-1}=R$). Hence, as $a_0=a_{n-1}=0$, we have:

$$\ell(\mathcal{S})=\sum_{i=1}^{n-1}|a_{i-1}-a_i|+ \sum_{j=0}^{n-1} |b_{j}|\leq \sum_{i=1}^{n-1}(|a_{i-1}|+|a_i|) +\sum_{j=0}^{n-1} |b_{j}|= \sum_{i=0}^{n-1} (2|a_i| + |b_i|).$$

\hfill $\blacksquare$

\begin{corollary}\label{corcost} Let $\mathcal{S}$ be a curve system on $D_n$ represented by the reduced Dynnikov coordinates $(a_0,b_0,\dots, a_{n-1},b_{n-1})$. Then, the length of the minimal standardizer of $\mathcal{S}$ is at most $$\dfrac{1 }{2}\sum_{i=0}^{n-1} (2|a_i| + |b_i|)(n-1)^2.$$
\end{corollary}

\emph{Proof.} By Lemma~\ref{s}, the length of the minimal standardizer of $\mathcal{S} $ is at most $\frac{1}{2}\ell(\mathcal{S})(n-1)^2$. Consider the bound for $\ell(\mathcal{S})$ given in Proposition~\ref{cotaminimalstand} and the result will follow. \hfill $\blacksquare$

\medskip

\begin{figure}[ht]
  \centering
  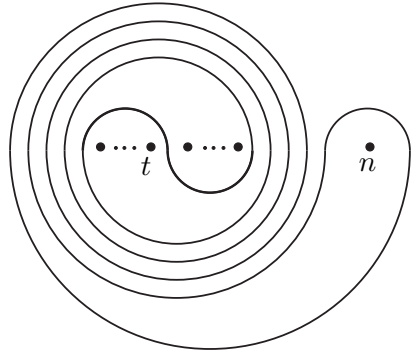
  \caption{The curve $\mathcal{S}_5$. }
  \label{optimalbound}
\end{figure}

\begin{remark}
To check that this bound is computationally optimal we need to find a case where at each step we can only remove a single bending point, i.e., we want to find a family of curve systems $\{\mathcal{S}_k\}_{k>0}$ such that the length of the minimal standardizer of $\mathcal{S}_k$ is quadratic on $n$ and linear on $\ell(\mathcal{S})$. Let $n=2t+1,\, t\in \mathbb{N}$. Consider the following curve system on $D_n$,
$$\mathcal{S}_0=\left\{t \smile n \frown\right\}$$
and the braid $\alpha=s(0,t,n-1)$. Now define $\mathcal{S}_k=(\mathcal{S}_0)^{\alpha^{-k}}$. The curve $\mathcal{S}_k$ is called a spiral with $k$ twists (see Figure~\ref{optimalbound}) and is such that $\ell(\mathcal{S}_k)= 2(k+1)$. Using Algorithm~\ref{algo}, we obtain that the minimal standardizer of this curve is $\alpha^k$, which has $k\cdot t ^2$ factors. Therefore, the number of factors of the minimal standardizer of $\mathcal{S}_k$ is of order $O(\ell(\mathcal{S}_k)\cdot n^2)$. 

\end{remark}

\begin{corollary}\label{bound}
Let $\mathcal{S}$ be a curve system on $D_n$ represented by the reduced Dynnikov coordinates $(a_0,b_0,\dots, a_{n-1},b_{n-1})$. Let $m=\sum_{i=0}^{n-1}(|a_i|+|b_i|)$. Then, the complexity of computing the minimal standardizer of $\mathcal{S}$ is $O(n^2 m \log(m))$.
\end{corollary}

\emph{Proof.}
First notice that $$\ell(\mathcal{S})\leq \sum_{i=0}^{n-1} (2|a_i| + |b_i|) \leq 2 \sum_{i=0}^{n-1}(|a_i|+|b_i|)=2m,$$ and that the transformation described in Proposition~\ref{transformation} involves a finite number of basic operations (addition and max). Applying $\sigma_j$ to the Dynnikov coordinates modifies only four such coordinates, and each maximum or addition between two numbers is linear on the number of digits of its arguments. This means that applying $\sigma_j$ to the curve has a cost of $O(\log(M))$, where $M=\max\{|a_i|, |b_i|\,|\, i=0,\dots , n-1\}$. By Corollary~\ref{corcost}, the number of iterations performed by the algorithm is $O(n^2m)$. Hence, as $M\leq m$, computing the minimal standardizer of $\mathcal{S}$ has complexity $O(n^2 m \log(m))$.  \hfill $\blacksquare$

%Thus, the algorithm is exponential on the size of $c$, which is $$\ell(c)= \sum_{i=0}^{n-1} (\log|a_i|+\log|b_i|)\leq 2(n-1)\log(M)\leq 2(n-1)\log(m).$$

\bigskip
To find the complexity of the algorithm which computes the minimal standardizer of a parabolic subgroup $P=(X, \alpha)$ of an Artin-Tits group $A$, we only need to know the cost of computing the $pn$-normal form of~$c_{P}$. If $x_r\cdots x_1 \Delta^{-p}$ with $p>0$ is the right normal form of~$c_{P}$, then its $pn$-normal form  is $(x_r\cdots x_{p+1})(x_p\cdots x_1 \Delta^{-p})$. Hence, we just have to compute the right normal form of~$c_{P}$ in order to compute the minimal standardizer. It is well known that this computation has quadratic complexity (for a proof, see \citep{Dehornoy2008a}). Thus, we have the following:

\begin{proposition}
Let $P=(X, \alpha)$ be a parabolic subgroup of an Artin-Tits group of spherical type, and let $\ell= \ell(\alpha)$ be the canonical length $\alpha$. Computing the minimal standardizer of $P$ has a cost of $O(\ell^2)$.

\end{proposition}

\bigskip

{\bf Acknowledgements. } This research was supported by a PhD contract founded by Universit\'{e} Rennes 1, Spanish Projects MTM2013-44233-P, MTM2016-76453-C2-1-P, FEDER and French-Spanish Mobility programme ``M\'{e}rim\'{e}e 2015". I also thank my PhD advisors, Juan Gonz\'alez-Meneses and Bert Wiest, for helping and guiding me during this research work.

\bibliography{Minimal}

\bigskip\bigskip{\footnotesize%
 \textit{ Mar\'{i}a Cumplido, UFR Math\'ematiques, Universit\'e de Rennes 1, France, and Departamento de \'Algebra, Universidad de Sevilla, Spain} \par
  \textit{E-mail address:} \texttt{\href{mailto:maria.cumplidocabello@univ-rennes1.fr}{maria.cumplidocabello@univ-rennes1.fr}}

}

\end{document}